\documentstyle[12pt]{article}
\begin{document}
\newcommand{\text}[1]{\mbox{{\rm #1}}}
\newcommand{\gd}{\delta}
\newcommand{\itms}[1]{\item[[#1]]}
\newcommand{\nin}{\in\!\!\!\!\!/}
\newcommand{\sub}{\subset}
\newcommand{\cntd}{\subseteq}
\newcommand{\go}{\omega}
\newcommand{\Pa}{P_{a^\nu,1}(U)}
\newcommand{\fx}{f(x)}
\newcommand{\fy}{f(y)}
\newcommand{\gD}{\Delta}
\newcommand{\gl}{\lambda}
\newcommand{\gL}{\Lambda}
\newcommand{\half}{\frac{1}{2}}
\newcommand{\sto}[1]{#1^{(1)}}
\newcommand{\stt}[1]{#1^{(2)}}
\newcommand{\Z}{\hbox{\sf Z\kern-0.720em\hbox{ Z}}}
\newcommand{\singcolb}[2]{\left(\begin{array}{c}#1\\#2
\end{array}\right)}
\newcommand{\ga}{\alpha}
\newcommand{\gb}{\beta}
\newcommand{\gga}{\gamma}
\newcommand{\ul}{\underline}
\newcommand{\ol}{\overline}
\newcommand{\qed}{\kern 5pt\vrule height8pt width6.5pt depth2pt}
\newcommand{\Lrraro}{\Longrightarrow}
\newcommand{\Nb}{|\!\!/}
\newcommand{\NN}{{\rm I\!N}}
\newcommand{\bsl}{\backslash}
\newcommand{\gt}{\theta}
\newcommand{\op}{\oplus}
\newcommand{\C}{{\bf C}}
\newcommand{\Q}{{\bf Q}}
\newcommand{\Op}{\bigoplus}
\newcommand{\CR}{{\cal R}}
\newcommand{\tr}{\bigtriangleup}
\newcommand{\grr}{\omega_1}
\newcommand{\ben}{\begin{enumerate}}
\newcommand{\een}{\end{enumerate}}
\newcommand{\ndiv}{\not\mid}
\newcommand{\bab}{\bowtie}
\newcommand{\hal}{\leftharpoonup}
\newcommand{\har}{\rightharpoonup}
\newcommand{\ot}{\otimes}
\newcommand{\OT}{\bigotimes}
\newcommand{\bwe}{\bigwedge}
\newcommand{\gep}{\varepsilon}
\newcommand{\gs}{\sigma}
\newcommand{\rbraces}[1]{\left( #1 \right)}
\newcommand{\bbox}{$\;\;\rule{2mm}{2mm}$}
\newcommand{\sbraces}[1]{\left[ #1 \right]}
\newcommand{\bbraces}[1]{\left\{ #1 \right\}}
\newcommand{\OO}{_{(1)}}
\newcommand{\TT}{_{(2)}}
\newcommand{\FF}{_{(3)}}
\newcommand{\minus}{^{-1}}
\newcommand{\CV}{\cal V}
\newcommand{\CVs}{\cal{V}_s}
\newcommand{\un}{U_q(sl_n)'}
\newcommand{\on}{O_q(SL_n)'}
\newcommand{\slq}{U_q(sl_2)}
\newcommand{\olq}{O_q(SL_2)}
\newcommand{\UU}{U_{(N,\nu,\go)}}
\newcommand{\HH}{H_{n,q,N,\nu}}
\newcommand{\ct}{\centerline}
\newcommand{\bs}{\bigskip}
\newcommand{\ms}{\medskip}
\newcommand{\noin}{\noindent}
\newcommand{\mat}[1]{$\;{#1}\;$}
\newcommand{\raro}{\rightarrow}
\newcommand{\map}[3]{{#1}\::\:{#2}\raro{#3}}
\newcommand{\alg}{{\rm Alg}}
\def\newtheorems{\newtheorem{theorem}{Theorem}[subsection]
                 \newtheorem{cor}[theorem]{Corollary}
                 \newtheorem{prop}[theorem]{Proposition}
                 \newtheorem{lemma}[theorem]{Lemma}
                 \newtheorem{defn}[theorem]{Definition}
                 \newtheorem{Theorem}{Theorem}[section]
                 \newtheorem{Corollary}[Theorem]{Corollary}
                 \newtheorem{Proposition}[Theorem]{Proposition}
                 \newtheorem{Lemma}[Theorem]{Lemma}
                 \newtheorem{Defn}[Theorem]{Definition}
                 \newtheorem{Definition}[Theorem]{Definition}
                 \newtheorem{Example}[Theorem]{Example}
                 \newtheorem{Remark}[Theorem]{Remark}
                 \newtheorem{claim}[theorem]{Claim}
                 \newtheorem{sublemma}[theorem]{Sublemma}
                 \newtheorem{example}[theorem]{Example}
                 \newtheorem{remark}[theorem]{Remark}
                 \newtheorem{question}[theorem]{Question}
                 \newtheorem{Question}[Theorem]{Question}
                 \newtheorem{conjecture}{Conjecture}[subsection]}
\newtheorems
\newcommand{\proof}{\par\noindent{\bf Proof:}\quad}
\newcommand{\dmatr}[2]{\left(\begin{array}{c}{#1}\\
                            {#2}\end{array}\right)}
\newcommand{\doubcolb}[4]{\left(\begin{array}{cc}#1&#2\\
#3&#4\end{array}\right)}
\newcommand{\qmatrl}[4]{\left(\begin{array}{ll}{#1}&{#2}\\
                            {#3}&{#4}\end{array}\right)}
\newcommand{\qmatrc}[4]{\left(\begin{array}{cc}{#1}&{#2}\\
                            {#3}&{#4}\end{array}\right)}
\newcommand{\qmatrr}[4]{\left(\begin{array}{rr}{#1}&{#2}\\
                            {#3}&{#4}\end{array}\right)}
\newcommand{\smatr}[2]{\left(\begin{array}{c}{#1}\\
                            \vdots\\{#2}\end{array}\right)}

\newcommand{\ddet}[2]{\left[\begin{array}{c}{#1}\\
                           {#2}\end{array}\right]}
\newcommand{\qdetl}[4]{\left[\begin{array}{ll}{#1}&{#2}\\
                           {#3}&{#4}\end{array}\right]}
\newcommand{\qdetc}[4]{\left[\begin{array}{cc}{#1}&{#2}\\
                           {#3}&{#4}\end{array}\right]}
\newcommand{\qdetr}[4]{\left[\begin{array}{rr}{#1}&{#2}\\
                           {#3}&{#4}\end{array}\right]}

\newcommand{\qbracl}[4]{\left\{\begin{array}{ll}{#1}&{#2}\\
                           {#3}&{#4}\end{array}\right.}
\newcommand{\qbracr}[4]{\left.\begin{array}{ll}{#1}&{#2}\\
                           {#3}&{#4}\end{array}\right\}}

\title{Some Properties and Examples of Triangular Pointed Hopf Algebras}
\author{
Shlomo Gelaki\\
Technion-Israel Institute of Technology\\
Department of Mathematics\\
Haifa, Israel 32000}
\maketitle

\section{Introduction}
A fundamental problem in the theory of Hopf algebras is the classification
and construction of finite-dimensional (minimal) triangular Hopf
algebras $(A,R)$ introduced by
Drinfeld [D]. Only recently this problem was completely solved for
semisimple $A$ over algebraically closed fields 
of characteristics $0$ and $p>>dim(A)$ (any $p$
if one assumes
that $A$ is also cosemisimple) [EG2,EG3].

In this paper we take the first step towards solving this problem for
finite-dimensional pointed Hopf algebras over an algebraically closed 
field $k$ of characteristic $0.$ 

We first prove that the fourth power of the antipode of any triangular
pointed Hopf algebra $A$ is the identity.
We do that by focusing on {\em minimal} triangular pointed Hopf
algebras $(A,R)$
(every triangular Hopf algebra contains a
minimal triangular sub Hopf algebra) and proving that the group algebra of
the group of grouplike elements of $A$ (which must be abelian) admits a
minimal triangular structure and consequently that $A$ has the structure
of a biproduct [R1]. We also generalize our result on the order of the
antipode to any finite-dimensional quasitriangular Hopf algebra $A$
whose Drinfeld element $u$ acts as a scalar in any irreducible
representation of $A$ (e.g. when $A^*$ is pointed). 

Second, we describe a method of construction of finite-dimensional
pointed Hopf algebras which admit a minimal triangular structure, and
classify all their minimal triangular structures. 

We conclude the paper by proving that {\em any} minimal triangular
Hopf
algebra which is generated as an algebra by grouplike elements and skew
primitive
elements is isomorphic to a minimal triangular pointed Hopf algebra
constructed
using our method.

Throughout the paper, unless otherwise specified, the ground field $k$
will be assumed to be algebraically closed with characteristic $0.$

\noin
{\bf Acknowledgment:} I am grateful to Pavel Etingof for numerous
conversations, for reading the manuscript and for making suggestions which
improved the presentation of the paper. 
\section{Pointed Hopf Algebras} 
The Hopf algebras which are studied in this paper are pointed. 
Recall that a Hopf algebra $A$ is pointed if its simple subcoalgebras are 
all $1-$dimensional or equivalently (when $A$ is finite-dimensional) if 
the irreducible representations of $A^*$ are all $1-$dimensional. 
Let $G(A)$ denote the group of grouplike elements of $A.$
For any $g,h\in G(A),$ we denote the vector space of $g:h$ skew primitives
of $A$ by $P_{g,h}(A):=\{x\in A|\gD(x)=x\ot g+h\ot x\}.$
Thus the classical primitive elements of $A$ are $P(A):=P_{1,1}(A).$
The element $g-h$ is always $g:h$ skew primitive. Let
$P'_{g,h}(A)$ denote a complement of $sp_k\{g-h\}$ in  
$P_{g,h}(A)$. Taft-Wilson theorem [TW] states that  
the first term $A_1$ of the coradical filtration of $A$ is given by:
\begin{equation}\label{wilson}
A_1=kG(A)\bigoplus (\Op_{g,h\in G(A)} P'_{g,h}(A)).
\end{equation}
In particular, if $A$ is {\em not} cosemisimple then there
exists $g\in G(A)$ such that $P'_{1,g}(A)\ne 0.$ 

If $A$ is a Hopf algebra over the field $k,$ which is generated
as an algebra by a subset $S$ of $G(A)$ and by $g:g'$ skew primitive
elements, where $g,g'$ run over $S,$ then $A$ is pointed and $G(A)$ is
generated as a group by $S$ (see e.g. [R4, Lemma 1]). 
\section{The Antipode of Triangular Pointed Hopf Algebras} 
In this section
we recall some of the properties of finite-dimensional triangular Hopf
algebras, and prove some new general results about 
finite-dimensional triangular pointed Hopf algebras.

Let $(A,R)$ be a finite-dimensional triangular Hopf algebra over 
a field $k$ of characteristic $0.$ Recall that the $R-$matrix $R$ 
satisfies the relation $R^{-1}=R_{21}$ or equivalently, the Drinfeld 
element $u$ of $A$ is a grouplike element [D]. If $A$ is semisimple
then 
$u^2=1$ (see e.g. [EG1]). The associated map
$f_R:A^{*cop}\raro A$ defined by $f_R(p)=(p\ot I)(R)$ is a Hopf
algebra
homomorphism which satisfies
\begin{equation}\label{qt5}
\sum <p_{(1)},a_{(2)}>a_{(1)}f_R(p_{(2)})=\sum 
<p_{(2)},a_{(1)}>f_R(p_{(1)})a_{(2)},\;p\in A^*,\;a\in A
\end{equation}
and
\begin{equation}\label{t}
\sum f_R(p_{(1)})f_R^*(p_{(2)})=<p,1>1,\;p\in A^*.
\end{equation}
Observe that (\ref{qt5}) is equivalent to the condition
$\Delta^{cop}(a)R=R\Delta(a)$ for all $a\in A,$ and (\ref{t}) is
equivalent to the condition $RR_{21}=1\ot 1.$

Let $A_R\subseteq A$ be the minimal sub Hopf algebra of $A$ corresponding
to $R$ [R2]. It is straightforward to verify that the corresponding map
$f_R:A_R^{*cop}\raro A_R$ defined by 
$f_R(p)=(p\ot I)(R)$ is a Hopf algebra isomorphism. This 
property of minimal triangular Hopf algebras will play a central role in
this paper. It implies in particular that $G(A_R)\cong G(A_R^*),$
and hence that the group $G(A_R)$ is abelian (see e.g. [G2]).
Thus, $G(A_R)\cong G(A_R)^*$ (where $G(A_R)^*$ denotes the character group of
$G(A_R)$), and we can identify the Hopf algebras
$k[G(A_R)^*]$ and $k[G(A_R)]^*.$  
Also, if $(A,R)$ is minimal triangular and pointed then $f_R$
being an isomorphism implies that $A^*$ is pointed as well.

Note that if $(A,R)$ is (quasi)triangular and $\pi:A\raro A'$ is an onto
morphism 
of Hopf algebras, then $(A',R')$ is (quasi)triangular as well, 
where $R'=(\pi\ot \pi)(R).$ 

We are ready now to prove our first results.
\begin{Theorem}\label{gen}
Let $(A,R)$ be a minimal triangular pointed Hopf algebra with Drinfeld 
element $u$ over an
algebraically closed field $k$ of characteristic $0,$
and set $K:=k[G(A)].$ Then there exists a projection of Hopf
algebras $\pi:A\raro K,$ and consequently
$A=B\times K$ is a biproduct where $B=\{x\in A|(I\ot 
\pi)\Delta(x)=x\ot 1\}\subseteq A.$ Moreover, $K$ admits a minimal 
triangular structure with Drinfeld element $u_{_K}=u.$
\end{Theorem}
\proof We know that $k[G(A)]^*\cong k[G(A)]$ and $k[G(A)]\cong 
k[G(A^{*cop})]$ as Hopf algebras, and hence that  
$dim(k[G(A)]^*)=dim(k[G(A^{*})])$. Consider the series of Hopf
algebra homomorphisms
$k[G(A)]\stackrel{i}{\hookrightarrow}A^{cop}\stackrel{f_R^{-1}}
{\longrightarrow}A^*\stackrel{i^*}{\longrightarrow}k[G(A)]^*,$
where $i$ is the inclusion map.
Since $A^*$ is pointed it follows from [M, 5.3.5] and the above remarks 
that $i^*_{|k[G(A^{*})]}:k[G(A^{*})]\raro 
k[G(A)]^*$ is an isomorphism of Hopf algebras, and 
hence that $i^*\circ f_R^{-1}\circ i$ determines a minimal quasitriangular 
structure on $k[G(A)]^*.$ This structure is in fact triangular since
$f_R^{-1}$ determines a triangular structure on $A^*.$
Clearly, $(i^*\circ f_R^{-1}\circ 
i)(u)=u_{_{K^*}}^{-1}=u_{_{K^*}}$ is the Drinfeld element of $K^*.$
Since $k[G(A)]$ and $k[G(A)]^*$ are isomorphic as Hopf algebras we 
conclude that $K$ admits a minimal triangular structure with Drinfeld
element $u_{_K}=u.$

Finally, set $\varphi:=i^*\circ f_R^{-1}\circ i$ and
$\pi:=\varphi^{-1}\circ i^* 
\circ f_R^{-1}.$ Then $\pi:A\raro k[G(A)]$ is onto, and moreover 
$\pi\circ i=\varphi^{-1}\circ i^*\circ f_R^{-1}\circ i=\varphi^{-1}\circ 
\varphi=id_{k[G(A)]}.$ Hence $\pi$ is a projection of Hopf algebras and 
by [R1], $A=B\times K$ is a biproduct where $B=\{x\in A|(I\ot 
\pi)\Delta(x)=x\ot 1\}$ as desired. This concludes the proof of the 
theorem. \qed
\begin{Theorem}\label{ant}
Let $(A,R)$ be any triangular pointed Hopf algebra with antipode $S$
and Drinfeld element $u$ over
a field $k$ of characteristic $0.$  
Then $S^4=I$ is the identity map of $A.$ If in addition $A_R$ is 
not semisimple and $A$ is finite-dimensional then $dim(A)$ is divisible by
$4.$
\end{Theorem}
\proof We may assume that $k$ is algebraically closed.
By [D], $S^2(x)=uxu^{-1}$ for all 
$x\in A.$ Let $K:=k[G(A_R)].$ Since $u\in A_R,$ and by Theorem 
\ref{gen}, $u=u_{_K}$ and $u_{_K}^2=1,$ we have that
$S^4=I.$ 

In order to prove the second claim, we may assume that $(A,R)$ 
is minimal (since by [NZ], $dim(A_R)$ divides $dim(A)$).
Since $A$ is not semisimple it follows from [LR] that $S^2\ne I,$ 
and hence that $u\ne 1.$ In particular, $|G(A)|$ is even. Now, let $B$ be as
in Theorem \ref{gen}. Since $S^2(B)=B,$ $B$ has a basis
$\left \{a_i,b_j|S^2(a_i)=a_i,S^2(b_j)=-b_j,1\le i\le n,1\le j\le m\right
\}.$
Hence by Theorem \ref{gen}, 
$$\left \{a_ig,b_jg|g\in G(A),1\le i\le n,1\le j\le m\right \}$$ 
is a basis of $A.$ Since by [R3], $tr(S^2)=0,$ we have that 
$0=tr(S^2)=|G(A)|(n-m),$ which implies that $n=m,$ and hence that $dim(B)$
is even as well. \qed

In fact, the first part of Theorem \ref{ant} can be generalized.
\begin{Theorem}\label{p}
Let $(A,R)$ be a finite-dimensional quasitriangular Hopf algebra with
antipode $S$ over a field $k$ of characteristic $0,$ and suppose that
the Drinfeld element $u$ of $A$ acts as a scalar in any irreducible
representation of $A$ (e.g. when $A^*$ is pointed). Then $u=S(u)$ and in
particular $S^4=I.$
\end{Theorem}
\proof We may assume that $k$ is algebraically closed. In any irreducible
representation $V$ of $A,$ $tr_{|V}(u)=tr_{|V}(S(u))$ (see
[EG1]). Since $S(u)$ also acts as a scalar in $V$ (the dual of
$S(u)_{|V}$ equals $u_{|V^*}$) it follows that
$u=S(u)$ in any irreducible representation of $A.$ Therefore, there exists
a basis of $A$ in which the operators of left multiplication by 
$u$ and $S(u)$ are represented by upper triangular matrices with the same
main diagonal. Hence the special grouplike element $uS(u)^{-1}$ is
unipotent. Since it has a finite order we conclude that $uS(u)^{-1}=1,$
and hence that $S^4=I.$ \qed
\begin{Remark}\label{ge}
If $(A,R)$ is a minimal triangular pointed Hopf algebra then all its
irreducible representations are $1-$dimensional. Hence Theorem \ref{p} is
applicable, and the first part of Theorem \ref{ant} follows.
\end{Remark}
\begin{Example}\label{sweed}
Let $A$ be Sweedler's $4-$dimensional Hopf algebra. It is generated as
an algebra by a
grouplike element $g$ and a $1:g$ skew primitive element $x$ satisfying
the relations
$g^2=1,$ $x^2=0$ and $gx=-xg.$ It is known [R2] that $A$ admits minimal
triangular structures all of which with $g$ as the Drinfeld element. In
this example, $K=k[<g>]$ and $B=sp\{1,x\}.$ Note that $g$ is central in $K$
but is not central in $A,$ so $(S_{|K})^2=I_K$ but $S^2\ne I$ in $A.$
However, $S^4=I.$
\end{Example}

Theorem \ref{ant} motivates the following two questions.
\begin{Question}\label{q}
{\rm 1) Is the fourth power of the antipode of any finite-dimensional
triangular Hopf algebra equal to the
identity? 2) Does the Drinfeld element $u$ of any finite-dimensional
triangular Hopf algebra satisfy $u^2=1$?}
\end{Question}

Note that a positive answer to Question \ref{q} 2) will imply that an
odd-dimensional triangular Hopf algebra must be semisimple.
\section{Construction of Minimal Triangular Pointed Hopf Algebras}
In this section we give a method for the construction of minimal
triangular pointed Hopf algebras which are {\em not}
necessarily semisimple.

Let $G$ be a finite abelian group, and $F:G\times G\raro k^*$ be a 
non-degenerate skew symmetric bilinear form on $G.$ That is,
$F(xy,z)=F(x,z)F(y,z),$ $F(x,yz)=F(x,y)F(x,z),$ $F(1,x)=F(x,1)=1,$ 
$F(x,y)=F(y,x)^{-1}$ for all $x,y,z\in G,$ and the map
$f:G\raro G^*$ defined by $<f(x),y>=F(x,y)$ for
all $x,y\in G$ is an isomorphism. Let
$U_F:G\raro \{-1,1\}$ be defined by $U_F(g)=F(g,g).$ Then $U_F$ is
a homomorphism of groups. Denote $U_F^{-1}(-1)$ by $I_F.$

\begin{Definition}\label{datum}
Let $k$ be an algebraically closed field of characteristic zero.
A datum \linebreak
${\cal D}=(G,F,n)$ is a triple where
$G$ is a finite abelian group, $F:G\times G\raro k^*$ 
is a non-degenerate skew symmetric bilinear form on $G,$ and $n$ is a
non-negative integer function $I_F\raro \Z^+,$ $g\mapsto n_g.$
\end{Definition}

\begin{Remark}\label{form}
1) The map $f:k[G]\raro k[G^*]$ determined by $<f(g),h>=F(g,h)$ for 
all $g,$ $h\in G$ determines a minimal triangular structure on $k[G^*].$

\noin
2) If $I_F$ is not empty then $G$ has an even order.
\end{Remark}

To each datum ${\cal D}$ we associate a Hopf algebra $H({\cal D})$ in the
following way. For each $g\in I_F,$ 
let $V_g$ be a vector space of dimension $n_g,$ and let
${\cal B}=\bigoplus_{g\in I_F} V_g.$
Then $H({\cal D})$ is generated as an algebra by $G\cup {\cal B}$ with the
following additional relations (to those of the group $G$ and the vector
spaces $V_g$'s):
\begin{equation}\label{3}
xy=F(h,g)yx\;\;and\;\;xa=F(a,g)ax 
\end{equation}
for all $g,h\in I_F,$ $x\in V_g,$ $y\in V_h$ and $a\in G.$

The coalgebra structure of $H({\cal D})$ is determined by letting $a\in G$
be a grouplike element and $x\in V_g$ be a $1:g$ skew primitive element
for all $g\in I_F.$ In particular, $\varepsilon(a)=1$ and
$\varepsilon(x)=0$ for all $a\in G$ and $x\in V_g.$

In the special case where $G=\Z_2=\{1,g\},$ $F(g,g)=-1$ and $n:=n_g,$
the associated Hopf algebra will be denoted by $H(n).$ Clearly,
$H(0)=k\Z_2,$ $H(1)$ is Sweedler's $4-$dimensional Hopf algebra, and
$H(2)$ is the $8-$dimensional Hopf algebra studied in [G1, Section 2.2] in
connection with KRH invariants of knots and $3-$manifolds. We remark that the
Hopf algebras $H(n)$ are studied in [PO1,PO2] where they are denoted by
$E(n).$

For a finite-dimensional vector space $V$ we let $\bigwedge V$ denote the
exterior algebra of $V.$ Set $B:=\bigotimes_{g\in I_F}\bigwedge V_g.$
\begin{Proposition}\label{dual}
1) The Hopf algebra $H({\cal D})$ is pointed and $G(H({\cal D}))=G.$

\noin
2) $H({\cal D})=B\times k[G]$ is a biproduct. 

\noin
3) $H({\cal D})_1= k[G]\bigoplus (k[G]{\cal B}),$ and $P_{a,b}(H({\cal
D}))=sp\{a-b\}\bigoplus aV_{a^{-1}b}$ for all $a,b\in G$
(here we agree that $V_{a^{-1}b}=0$ if $a^{-1}b\notin I_F$).
\end{Proposition}
\proof 
Part 1) follows since (by definition) $H({\cal D})$ is generated as an
algebra by grouplike elements and skew primitive elements.
Now, it is straightforward to verify that the map
$\pi:H({\cal D})\raro k[G]$ determined by
$\pi(a)=a$ and $\pi(x)=0$ for all $a\in G$ and $x\in
{\cal B},$ is a projection of Hopf algebras. Since $B=\{x\in H({\cal
D})|(I\ot
\pi)\Delta(x)=x\ot 1\},$ Part 2) follows from [R1]. Finally,
by Part 2), $B$ is a braided graded Hopf algebra in the Yetter-Drinfeld
category $_{k[G]}^{k[G]}{\cal YD}$ (see e.g. [AS]) with respect to the
grading where the elements of ${\cal B}$ are homogeneous of
degree $1.$ Write $B=\bigoplus_{n\ge 0}B(n),$ where $B(n)$
denotes the homogeneous component of degree $n.$ Then, $B(0)=k1=B_1$
(since $B\cong H({\cal D})/H({\cal D})k[G]^+$ as coalgebras, it is
connected). Furthermore, by similar arguments used in the proof of
[AS, Lemma 3.4], $P(B)=B(1)={\cal B}.$ But
then by [AS, Lemma 2.5], $H({\cal D})$ is coradically graded (where the
$nth$ component $H({\cal D})(n)$ is just $B(n)\times k[G]$) which means
by definition that $H({\cal D})_1= H({\cal D})(0)\bigoplus H({\cal
D})(1)=k[G]\bigoplus (k[G]{\cal B})$ as desired. The second statement
of Part 3) follows now, using (1), by counting dimensions. \qed

In the following we determine {\em all} the minimal triangular structures
on $H({\cal D}).$ Let \linebreak $f:k[G]\raro k[G^*]$ be the isomorphism
from
Remark \ref{form} 1), and set $I_F':=\{g\in I_F|n_g\ne 0\}.$ Let
$\Phi$ be the set of all isomorphisms $\varphi:G^*\raro G$
satisfying $\varphi^*(\ga)=\varphi(\ga ^{-1})$ for all $\ga\in G^*,$ and
$(\varphi\circ f)(g)=g$ for all $g\in I_F'.$

Extend any $\ga \in G^*$ to an algebra homomorphism 
$H({\cal D})\raro k$ by setting $\ga (z)=0$ for all $z\in {\cal B}.$
Extend any $x\in V_g^*$ to $P_x\in H({\cal D})^*$ by setting $<P_x,ay>=0$
for all $a\in G$ and $y\in \bigotimes_{g\in I_F}\bigwedge V_g$ of degree
different from $1,$ and $<P_x,ay>=\delta_{g,h}<x,y>$ for all $a\in G$ and
$y\in V_h.$ We shall identify the vector
spaces $V_g^*$ and $\{P_x|x\in V_g^*\}$ via the map $x\mapsto P_x.$

For $g\in I_F',$ let $S_g(k)$ be the set of all isomorphisms $M_g:V_g^*\raro
V_{g^{-1}}.$ Let $S(k)\subseteq \times _{g\in I_F'} S_g(k)$ be the set of all 
tuples $(M_g)$ satisfying $M_g^*=M_{g^{-1}}$ for all $g\in I_F'.$

\begin{Theorem}\label{quas}
1) For each $T:=(\varphi,(M_g))\in \Phi\times S(k),$ there exists a unique Hopf
algebra isomorphism $f_T:H({\cal
D})^{*cop}\raro H({\cal D})$ determined by $\ga \mapsto
\varphi(\ga)$ and $P_x\mapsto M_{g}(x)$ for $\ga\in G^*$ and $x\in
V_g^*.$  

\noin
2) There is a one to one correspondence between
$\Phi\times S(k)$ and the
set of minimal triangular structures on $H({\cal D})$ given by $T\mapsto
f_T.$
\end{Theorem}
\proof We first show that $f_T$ is a well defined isomorphism of Hopf
algebras.
Using Proposition \ref{dual} 2), it is straightforward to
verify that $\Delta(P_x)=\varepsilon\ot 
P_x+P_x\ot f(g^{-1}),$ $P_x\ga=<\ga,g>\ga P_x,$ and
$P_xP_y=F(h,g)P_yP_x$ for all $\ga\in G^*,$  
$g,h\in I_F$, $x\in V_g^*$ and $y\in V_h^*.$
Let ${\cal B}^*:=\{P_x|x\in V_g^*,\;g\in I_F\},$ and $H$ be the
sub Hopf algebra of $H({\cal D})^{*cop}$ generated as an algebra by
$G^*\bigcup {\cal B}^*.$
Then, using (4) and our assumptions on $T,$
it is straightforward to verify that the map $f_T^{-1}:H({\cal
D})\raro H$ determined by $a\mapsto \varphi^{-1}(a)$ and
$z\mapsto M_{g}^{-1}(z)$ for $a\in G$ and $z\in
V_{g^{-1}},$ is a surjective homomorphism of Hopf algebras.
Let us verify for instance
that $f_T^{-1}(za)=F(a,g)f_T^{-1}(az).$ Indeed,
this is equivalent to $<\varphi^{-1}(a),g>=<f(a),g>$ which in turn
holds by our assumptions on $\varphi.$ Now, using Proposition 4.3 3), it
is
straightforward to verify that $f_T^{-1}$ is injective on
$P_{a,b}(H({\cal D}))$ for all $a,b\in G.$
Since $H({\cal D})$ is pointed, $f_T^{-1}$ is also injective
(see e.g. [M, Corollary 5.4.7]). This implies that $H=H({\cal
D})^{*cop},$ and that $f_T:H({\cal D})^{*cop}\raro H({\cal D})$ is an
isomorphism of Hopf algebras as desired. Note that in particular,
$G^*=G(H({\cal D})^*).$

The fact that $f_T$ satisfies (\ref{qt5}) follows from 
a straightforward computation (using (\ref{3})) since it is enough to
verify it for algebra generators $p\in G^*\cup {\cal B}^*$ of $H({\cal
D})^{*cop},$ and $a\in G\cup {\cal B}$ of $H({\cal D}).$ 

We have to show that $f_T$ satisfies (\ref{t}).
Indeed, it is straightforward to verify that \linebreak $f_T^*:H({\cal
D})^{*op}\raro H({\cal D})$ is determined by $\ga \mapsto
\varphi(\ga^{-1})$ and $P_x\mapsto gM_{g}(x)$ for
$\ga\in 
G^*$ and $x\in V_g^*.$ Hence,
$f_T^*=f_T\circ S,$
where $S$ is the antipode of $H({\cal D})^*,$ as desired. 

We now have to show that {\em any} minimal triangular structure on $H({\cal
D})$ comes from $f_T$ for some $T.$ 
Indeed, let ${\bf f}:H({\cal D})^{*cop}\raro H({\cal D})$ 
be any Hopf isomorphism. Then ${\bf f}$ must map
$G^*$ onto $G,$ $\{f(g^{-1})|g\in I_F'\}$ onto $I_F',$ and 
$P_{f(g^{-1}),\varepsilon}(H({\cal D})^{*cop})$ bijectively onto
$P_{1,\varphi(f(g^{-1}))}(H({\cal D})).$ Therefore
there exists an invertible operator $M_g:V_{g}^*\raro
V_{\varphi(f(g^{-1}))}$ such that ${\bf f}$ is determined by $\ga\mapsto
\varphi(\ga)$ and $P_x\mapsto M_g(x).$ Suppose
${\bf f}$ satisfies (\ref{qt5}). Then letting $p=P_x$ and $a\in G$ in
(\ref{qt5}) yields that $a{\bf f}(P_x)=F(a,g){\bf f}(P_x)a$ for all $a\in
G.$ But by (\ref{3}), this is equivalent to $(\varphi\circ f)(g)=g$
for all
$g\in I_F'.$ Since by Theorem \ref{gen}, $\varphi:k[G^*]\raro k[G]$
determines a minimal triangular structure on $k[G]$ it follows that
$\varphi\in \Phi.$ Since ${\bf f}:H({\cal D})^{*cop}\raro H({\cal D})$ 
satisfies (\ref{t}), $(M_g)\in S(k),$
and hence ${\bf f}$ is of the form $f_T$ for some $T$ as desired. \qed 

For a triangular structure on $H({\cal D})$ corresponding to the map 
$f_T,$ we let $R_T$ denote the corresponding $R-$matrix. 

\begin{Remark}
{\rm Note that if $n_{g^{-1}}\ne n_g$ for some $g\in I_F',$ 
then $S(k)$ is empty and $H({\cal D})$ does {\em not} have
a minimal triangular structure.
}
\end{Remark}

\section{Classification of Minimal Triangular Hopf Algebras Generated as 
Algebras by Grouplike and Skew Primitive Elements}
In this section we use Theorems \ref{gen} and \ref{ant} to classify
minimal triangular Hopf algebras which are generated as algebras by 
grouplike elements and skew primitive elements. Namely, we prove:

\begin{Theorem}\label{main}
Let $(A,R)$ be a minimal triangular Hopf
algebra over an algebraically closed field $k$ of characteristic
$0.$ If $A$ is generated as an
algebra by grouplike elements and skew primitive elements, 
then there exist a datum ${\cal D}$ and $T\in \Phi\times S(k)$ such that
$(A,R)\cong (H({\cal D}),R_T)$ as triangular Hopf algebras.
\end{Theorem}

Before we prove Theorem \ref{main} we need to fix some notation and prove
a few lemmas.

In what follows, $(A,R)$ will always be a minimal triangular pointed Hopf
algebra over $k,$ $G:=G(A)$ 
and $K:=k[G(A)].$ For any $g\in G,$ $P_{1,g}(A)$ is a $<g>-$module 
under conjugation, and $sp\{1-g\}$ is a submodule of
$P_{1,g}(A).$ Let $V_g\subset P_{1,g}(A)$ be its complement, and set
$n_g:=dim(V_g).$

By Theorem \ref{gen}, $A=B\times K$ where $B=\{x\in A|(I\ot 
\pi)\Delta(x)=x\ot 1\}\subseteq A$ is a 
left coideal subalgebra of $A$ (equivalently, $B$ is an object in the
Yetter-Drinfeld category $_{k[G]}^{k[G]}{\cal Y}{\cal D}$).
Note that $B\cap K=k1.$ Let
$\rho:B\raro K\ot B$ be the associated comodule structure and write
$\rho(x)=\sum x^1\ot x^2,$ $x\in B.$ 
By [R1], $B\cong A/AK^+$ as coalgebras, hence $B$ is a connected pointed 
coalgebra. Let $P(B):=\{x\in B|\Delta_B(x)=x\ot 1+1\ot x\}$ be the space 
of primitive elements of $B.$
\begin{Lemma} \label{shtut1} 
For any $g\in G,$ $V_g=\{x\in P(B)|\rho(x)=g\ot x\}.$
\end{Lemma}
\proof 
Let $x\in V_g.$ Since $g$ acts on $V_g$ by 
conjugation we may assume by [G1, Lemma 0.2], that $gx=\go xg$ for some
$1\ne \go\in k.$ Since $\pi(x)$ and $\pi(g)=g$ commute we 
must have that $\pi(x)=0.$ But then $(I\ot \pi)\Delta(x)=x\ot 1$ and 
hence $x\in B.$
Since $\Delta(x)=\sum x_1\times x_2^1\ot x_2^2\times 1,$ applying the
maps $\varepsilon\ot I\ot I\ot \varepsilon$ and 
$I\ot \varepsilon\ot I\ot \varepsilon$ to both sides of the equation 
$\sum x_1\times x_2^1\ot x_2^2\times 1=x\times 1\ot 1\times 1+1\times g\ot 
x\times 1,$ yields that $x\in P(B)$ and $\rho(x)=g\ot x$ as desired.

Suppose that $x\in P(B)$ satisfies $\rho(x)=g\ot x.$ Since
$\Delta(x)=x\ot 1+\rho(x),$ it follows that $x\in V_g$ as desired. \qed

\begin{Lemma} \label{shtut2}
For all $x\in V_g,$ $x^2=0$ and $gx=-xg.$
\end{Lemma}
\proof Suppose $V_g\ne 0$ and
let $0\ne x\in V_g.$ Then $S^2(x)=g^{-1}xg,$
$g^{-1}xg\ne x$ by [G1, Lemma 0.2], and $g^{-1}xg\in V_g.$ Since by 
Theorem \ref{ant}, $S^4=I$ it follows that $g^2$ and $x$ commute, and
hence that $gx=-xg$ for all $x\in V_g.$

Second we wish to show that $x^2=0.$
By Lemma \ref{shtut1}, $x\in B$ and hence $x^2\in B$ ($B$ is a subalgebra
of $A$).
Since $\Delta(x^2)=x^2\ot 1+ g^2\ot x^2,$ and $x^2$ and
$g^2$ commute, it follows from [G1, Lemma 0.2] that 
$x^2=\ga(1-g^2)\in K$ for some $\ga\in k.$ We thus conclude that 
$x^2=0,$ as desired. \qed

Recall that the map $f_R:A^{*cop}\raro A$ is an isomorphism of Hopf
algebras, and let \linebreak
$F:G\times G\raro k^*$ be the associated non-degenerate skew
symmetric bilinear form on $G$ defined by
$F(g,h):=<f_R^{-1}(g),h>$ for all $g,h\in G.$

\begin{Lemma} \label{shtut3}
For any $x\in V_g$ and $y\in V_h,$ $xy=F(h,g)yx.$  
\end{Lemma}
\proof
If either $V_g=0$ or $V_h=0,$ there is nothing to prove. Suppose
$V_g,V_h\ne 0,$ and let $0\ne x\in V_g$ and $0\ne y\in V_h.$
Set $P:=f_R^{-1}(x).$ Then $P\in
P_{f_R^{-1}(g),\varepsilon}(A^{*cop}).$
Substituting $p:=P$ and $a:=y$ in equation (\ref{qt5}) yields that
$yx-F(g,h)xy=<P,y>(1-gh).$ Since $<P,y>(1-gh)\in B\cap K,$
it is equal to $0,$ and hence $yx=F(g,h)xy.$ \qed

\begin{Lemma} \label{shtut4}
For any $a\in G$ and $x\in V_g,$ $xa=F(a,g)ax.$
\end{Lemma} 
\proof Set $P:=f_R^{-1}(x).$ Then the result follows by  
letting $p:=P$ and $a\in G$ in (\ref{qt5}), and noting that $P\in
P_{f_R^{-1}(g), \varepsilon}(A^{*cop}).$ \qed

We can now prove Theorem \ref{main}.

\noin
{\bf Proof of Theorem \ref{main}:}
Let $n:I_F\raro \Z^+$ be the non-negative integer function defined by
$n(g)=n_g,$ and let ${\cal D}:=(G,F,n).$ Since by assumption, $A$ is
generated as an algebra by $G\cup (\bigoplus_{g\in I_F}V_g)$ it is pointed.
By Lemmas \ref{shtut2}-\ref{shtut4}, relations (\ref{3}) are satisfied.
Therefore there exists a surjection of Hopf algebras $\varphi:H({\cal
D})\raro A.$ Using Proposition 4.3 3), it is straightforward to verify
that $\varphi$ is injective on $P_{a,b}(H({\cal D}))$ for all
$a,b\in G.$ Since $H({\cal D})$ is pointed, $\varphi$ is also injective
(see e.g. [M, Corollary 5.4.7]). Hence $\varphi$ is an isomorphism of
Hopf algebras. The rest of the theorem follows now from
Theorem \ref{quas}. \qed

Theorem \ref{main} raises the natural question whether any minimal
triangular pointed Hopf algebra in characteristic $0$ is generated as an
algebra by grouplike elements and skew primitive elements. Recently, it
was shown that the answer to this question is positive [AEG]. Hence, Theorem
\ref{main} gives an explicit classification of minimal triangular Hopf algebras.


\begin{thebibliography}{[AEG] } 
\bibitem
[AEG]{aeg} N. Andruskiewitsch, P. Etingof and S. Gelaki,
Triangular Hopf 
Algebras with the Chevalley Property, {\em submitted},
math.QA/0008232.
\bibitem
[AS]{as} N. Andruskiewitsch and H.-J. Schneider, Lifting of Quantum
Linear Spaces and Pointed Hopf Algebras of Order $p^3$, {\em J.
Algebra} {\bf 209} (1998), 658-691.
\bibitem
[D]{d} V. Drinfeld, On Almost Cocommutative Hopf Algebras, {\em
Leningrad Mathematics Journal} {\bf 1} (1990), 321-342.
\bibitem
[EG1]{eg1} P. Etingof and S. Gelaki, Some Properties of Finite-Dimensional
Semisimple Hopf Algebras, {\em Mathematical Research Letters} {\bf 5}
(1998), 191-197.
\bibitem
[EG2]{eg2} P. Etingof and S. Gelaki, A Method of Construction of
Finite-Dimensional Triangular Semisimple Hopf Algebras, {\em Mathematical
Research Letters} {\bf 5} (1998), 551-561.
\bibitem
[EG3]{eg3} P. Etingof and S. Gelaki, The Classification of Triangular  
Semisimple and Cosemisimple Hopf Algebras Over an Algebraically Closed 
Field, {\em International Mathematics Research Notices} {\bf 5} (2000),
223-234. 
\bibitem [G1]{g1} S. Gelaki, On Pointed Ribbon Hopf Algebras, {\em Journal 
of Algebra} {\bf 181} (1996), 760-786.
\bibitem [G2]{g2} S. Gelaki, Quantum Groups of
Dimension $pq^2,$ {\em Israel Journal of Mathematics} {\bf 102} (1997),
227-267.
\bibitem [LR]{lr} R. Larson and D. Radford, Semisimple 
Cosemisimple Hopf Algebras, {\em American Journal of Mathematics} 
{\bf 110} (1988), 187-195.
\bibitem[M]{m} S. Montgomery, Hopf algebras and their actions on rings,
{\em CBMS Lecture Notes} {\bf 82}, AMS, 1993.
\bibitem [NZ]{nz} W. D. Nichols and M. B. Zoeller, A Hopf algebra freeness
theorem, {\em American Journal of Mathematics} {\bf 111} (1989), 381-385.
\bibitem [PO1]{po1} F. Panaite and F. V. Oystaeyen, Quasitriangular
structures for some pointed Hopf algebras of dimension $2^n,$ {\em
Communications in Algebra}, to appear.
\bibitem [PO2]{po2} F. Panaite and F. V. Oystaeyen, Clifford-type algebras as
cleft extensions for some pointed Hopf algebras, {\em Communications in
Algebra}, to appear. 
\bibitem [R1]{r1} D. E. Radford, The Structure of Hopf Algebras with a 
Projection, {\em Journal of Algebra} {\bf 2} (1985), 322-347.
\bibitem [R2]{r2} D. E. Radford, Minimal quasitriangular Hopf algebras, 
{\em Journal of Algebra} {\bf 157} (1993), 285-315.
\bibitem [R3]{r3} D. E. Radford, The trace function and Hopf
Algebras {\em Journal of Algebra} {\bf 163} (1994), 583-622.
\bibitem [R4]{r4} D. E. Radford, On Kauffman's Knot Invariants Arising
from
Finite-Dimensional Hopf algebras, Advances in Hopf algebras, 158, 
205-266, lectures notes in pure and applied mathematics, Marcel Dekker, 
N. Y., 1994.
\bibitem [TW]{tw} E. J. Taft and R. L. Wilson, On antipodes in pointed Hopf 
algebras, {\em Journal of Algebra} {\bf 29} (1974), 27-32.
\end{thebibliography}
\end{document}